\documentclass[11pt, draft]{amsart}
\usepackage{amssymb}
\usepackage[mathscr]{eucal}

\theoremstyle{plain}
\newtheorem*{theorem}{Theorem}

\theoremstyle{remark}
\newtheorem*{remark}{Remark}

\def\A{\mathcal A}
\def\B{\mathcal B}
\def\C{\mathbb C}
\def\la{\langle}
\def\lam{\lambda}
\def\N{\mathbb N}
\def\ot{\otimes}

\def\ra{\rangle}
\def\rank{\operatorname{rank}}

\begin{document}
\title[]{On isomorphisms of standard operator algebras}
\author{LAJOS MOLN\'AR}
\address{Institute of Mathematics and Informatics\\
         University of Debrecen\\
         Faculty of Natural Sciences\\
         4010 Debrecen, P.O.Box 12, Hungary}
\email{molnarl@math.klte.hu}
\thanks{  This research was supported by the
          Hungarian National Foundation for Scientific Research
          (OTKA), Grant No. T030082, T031995, and by
          the Ministry of Education, Hungary, Reg.
          No. FKFP 0349/2000}
\subjclass{Primary: 47B49}
\keywords{}
%\date{\today}
\begin{abstract}
The aim of this paper is to show that between standard operator algebras
every bijective map with a certain multiplicativity property related to
the Jordan triple isomorphisms of associative rings is automatically
additive.
\end{abstract}
\maketitle

\section{Introduction}
It is a surprising result of Martindale \cite[Corollary]{Martindale}
that every multiplicative bijective map from a prime ring containing a
nontrivial idempotent onto an arbitrary ring is necessarily additive.
Therefore, one can say that the multiplicative structure of rings of
that kind completely determines their ring structure.
This result has been utilized by \v Semrl in \cite{Semrl} to describe
the form of the semigroup isomorphisms of standard operator algebras on
Banach spaces.
The aim of this paper is to generalize this result quite significantly.
Other results on the additivity of multiplicative maps
(in fact, *-semigroup homomorphisms) between operator
algebras can be found in \cite{Hak, Molnar}.

Besides additive and multiplicative maps (that is, ring
homomorphisms) between rings, sometimes one has to consider Jordan
homomorphisms. The Jordan structure of associative rings has been
studied by many people in ring theory.
Moreover, Jordan operator
algebras have serious applications in the mathematical foundations of
quantum mechanics.
If $\mathcal R, \mathcal R'$ are rings and $\phi:\mathcal R\to
\mathcal R'$ is a transformation, then it is called a Jordan
homomorphism in case
\[
\phi(A+B)=\phi(A)+\phi(B)
\]
and
\[
\phi(AB+BA)=\phi(A)\phi(B)+\phi(B)\phi(A)
\]
holds for every $A,B\in \mathcal R$.
Clearly, every ring homomorphism is a Jordan homomorphism and the same
is true for the ring antihomomorphisms (the transformation $\phi:
\mathcal R \to \mathcal
R'$ is called a ring antihomomorphism if $\phi$ is additive and
satisfies $\phi(AB)=\phi(B)\phi(A)$ for all $A,B \in \mathcal R$).
One can go a little further by even more weakening the
multiplicativity property of the transformations in question as follows.
It is easy to see that if the ring $\mathcal R'$ is 2-torsion free (this
means
that $2A=0$ implies $A=0$), then every Jordan homomorphism $\phi:
\mathcal R \to \mathcal R'$ is a so-called
Jordan triple homomorphism, that is, $\phi$ is an additive function
satisfying
\begin{equation}\label{E:triple11}
\phi(ABA)=\phi(A)\phi(B)\phi(A) \qquad (A,B \in \mathcal R)
\end{equation}
(see the Introduction in \cite{Bresar}).
The aim of this paper is to show that in the same situation as in
\cite{Semrl}, that is, in the case of standard operator algebras acting
on infinite dimensional Banach spaces every bijective map satisfying the
equation \eqref{E:triple11} is automatically linear or conjugate-linear
and continuous.

\section{The Result}

We begin with the notation and definitions that we shall use throughout.

All linear spaces are considered over the complex field.
Let $X$ be a Banach space. Denote by $B(X)$ and $F(X)$ the algebra of
all bounded linear operators on $X$ and the ideal of all finite rank
operators in $B(X)$, respectively. A subalgebra of $B(X)$ which contains
$F(X)$ is called a standard operator algebra on $X$.
For any $n\in \N$, we denote by $M_n(\C)$ the algebra of all $n\times n$
complex matrices and ${}^t$ stands for the transpose.

The dual space of $X$ is denoted $X^*$ and
$A^*$ stands for the Banach space adjoint of the bounded linear
operator $A$ on $X$.
If $x\in X$ and $f\in X^*$, then $x\ot f$ denotes the operator defined
by
\[
(x\ot f)(z)=f(z)x \qquad (z\in X).
\]
Similarly, if $H$ is a Hilbert space and $x,y\in H$, then $x\ot y$
denotes the operator defined by
\[
(x\ot y)(z)=\la z,y\ra x \qquad (z\in H).
\]

Two idempotents $P,Q\in B(X)$ are
said to be mutually orthogonal (in the algebraic sense) if $PQ=QP=0$.
One can introduce a partial ordering on the set of all idempotents in
$B(X)$ by defining $P\leq Q$ if and only if $PQ=QP=P$.
An element $R$ of $B(X)$ is called a tripotent if $R^3=R$.

Now, the result of the paper reads as follows.

\begin{theorem}\label{T:triple2}
Let $X,Y$ be complex Banach spaces, $\dim X\geq 3$, and let
$\A \subset B(X)$ and $\B \subset B(Y)$ be
standard operator algebras.
Suppose that $\phi:\A \to \B$ is a bijective
transformation satisfying
\begin{equation}\label{E:triple3}
\phi(ABA)=\phi(A)\phi(B)\phi(A) \qquad (A,B \in \A).
\end{equation}
If $X$ is infinite dimensional, then we have the following
possibilities:
\begin{itemize}
\item[(i)]
there exists an invertible bounded linear operator $T: X\to Y$
and $c \in \{ -1, 1\}$ such that
\[
\phi(A)=cTAT^{-1} \qquad (A\in \A);
\]
\item[(ii)]
there exists an invertible bounded conjugate-linear operator $T: X\to
Y$ and $c \in \{ -1, 1\}$ such that
\[
\phi(A)=cTAT^{-1} \qquad (A\in \A);
\]
\item[(iii)]
there exists an invertible bounded linear operator $T: X^*\to
Y$ and $c \in \{ -1, 1\}$ such that
\[
\phi(A)=cTA^*T^{-1} \qquad (A\in \A);
\]
\item[(iv)]
there exists an invertible bounded conjugate-linear operator $T: X^*\to
Y$ and $c \in \{ -1, 1\}$ such that
\[
\phi(A)=cTA^*T^{-1} \qquad (A\in \A).
\]
\end{itemize}
If $X$ is finite dimensional, then we have $\dim X=\dim Y$.
So, our bijective transformation $\phi$ can be supposed to act on a
matrix algebra $M_n(\C)$. In this case
we have the following possibilities:
\begin{itemize}
\item[(v)]
there exist a ring automorphism $h$ of $\C$, an invertible matrix
$T\in M_n(\C)$ and $c \in \{ -1, 1\}$ such that
\[
\phi(A)=cTh(A)T^{-1} \qquad (A\in M_n(\C));
\]
\item[(vi)]
there exist a ring automorphism $h$ of $\C$, an invertible matrix
$T\in M_n(\C)$ and $c \in \{ -1, 1\}$ such that
\[
\phi(A)=cTh(A)^t T^{-1} \qquad (A\in M_n(\C)).
\]
\end{itemize}
Here, $h(A)$ denotes the matrix obtained from $A$ by applying $h$ on
every entry of it.
\end{theorem}

\begin{remark}
According to the referee's wish we point out that there are a lot of
noncontinuous ring automorphisms of the complex field. See, for example,
\cite{Kestelman}.
\end{remark}

\begin{proof}[Proof of the theorem]
First note that $\phi$ preserves the tripotents in $\A$ and $\B$ in both
directions, that is, an operator $P\in \A$ is a tripotent if and only if
so is $\phi(P)$.

We show that for every $n=0, 1, \ldots$ and tripotent $P\in \A$ we have
$\rank P=n$ if and only if $\rank \phi(P)=n$.
First observe that $\phi(0)=0$. Indeed, since $i\phi(0)\in \B$, there
exists an $A\in \A$ such that $i\phi(0)=\phi(A)$. It follows that
$-\phi(0)=(i\phi(0))\phi(0)(i\phi(0))=\phi(A0A)=\phi(0)$ and this
implies that $\phi(0)=0$.
So, we have
the rank preserving property of $\phi$ for $n=0$.
It follows from the first part of the proof of
\cite[Theorem 4]{OmladicSemrl} that every tripotent on a Banach
space is the difference of two mutually orthogonal idempotents
(to be honest, the theorem in question is about Hilbert spaces, but
the part
of the proof that we need here applies also for Banach spaces).
Suppose now that $\phi(P)\in \B$ is a rank-$k$ tripotent
if and only if $P\in \A$ is a rank-$k$ tripotent holds true for
$k=0,\ldots,n$. Let $P\in \A$ be a rank-$(n+1)$ tripotent.
Then the rank of $\phi(P)$ is at least $n+1$.
Let $Q\in \B$ be a rank-$(n+1)$ tripotent such that
$\phi(P)Q\phi(P)=Q$ and $Q\phi(P)Q=Q$.
The existence of such a tripotent follows from the representation of
tripotents as differences of mutually orthogonal idempotents mentioned
above. Let $Q'=\phi^{-1}(Q)$.
We have
\begin{equation}\label{E:triple4}
PQ'P=Q' \quad \text{and} \quad Q'PQ'=Q'.
\end{equation}
Clearly, the rank of $Q'$ is at least $n+1$. On the other hand,
the first equality in \eqref{E:triple4} gives us that the range of $Q'$
is included in the range of $P$, so the rank of $Q'$ is exactly $n+1$.
The tripotent $P$ is the difference of two mutually orthogonal
idempotents. These idempotents induce a splitting of $X$ into the
direct sum of three closed subspaces.
With respect to this splitting every operator has a matrix representation.
In particular, we can write
\[
P=
\left[
\begin{matrix}
I & 0   & 0\\
0 & -I  & 0\\
0 & 0   & 0
\end{matrix}
\right].
\]
Let
\[
Q'=
\left[
\begin{matrix}
Q_{11} & Q_{12} & Q_{13}\\
Q_{21} & Q_{22} & Q_{23}\\
Q_{31} & Q_{32} & Q_{33}
\end{matrix}
\right]
\]
be the representation of $Q'$.
It follows from the first equality in \eqref{E:triple4} that
the only possibly nonzero entries in the matrix of $Q'$ are $Q_{11}$ and
$Q_{22}$. The second equality in \eqref{E:triple4} now implies that
$Q_{11}$ and $-Q_{22}$ are idempotents. By the equality of the ranks of
$P$ and $Q'$ we conclude that $P=Q'$.
Therefore, we have $\phi(P)=Q$ and this gives us that the rank of
$\phi(P)$ is $n+1$.
By symmetry, we obtain that if $\phi(P)$ has rank $n+1$, then the same
must be true for $P$.

The key step of the proof now follows. We intend to prove that
if $P',P''\in \A$ are mutually orthogonal rank-1 idempotents,
then $\phi(P'+P'')=\phi(P')+\phi(P'')$. In order to verify this, let
$P\in \A$ be a rank-3 idempotent. Then it follows that $Q=\phi(P)$ is a
rank-3 tripotent.
Let $Q=R_1-R_2$, where $R_1,R_2 \in B(Y)$ are idempotents with
$R_1R_2=R_2R_1=0$.
If $A\in \A$ is any operator satisfying $PAP=A$, then
we have $Q\phi(A)Q=\phi(A)$. We compute
\[
R_1\phi(A)R_2=
R_1(Q\phi(A)Q)R_2=
-R_1\phi(A)R_2
\]
which implies that $R_1\phi(A)R_2=0$. We similarly have
$R_2\phi(A)R_1=0$. It follows that
$\phi(A)=R_1\phi(A)R_1+R_2\phi(A)R_2$.
Conversely, if $\phi(A)=R_1\phi(A)R_1+R_2\phi(A)R_2$, then we
infer
\[
Q\phi(A)Q=
Q(R_1\phi(A)R_1+R_2\phi(A)R_2)Q=
\]
\[
R_1\phi(A)R_1+R_2\phi(A)R_2=
\phi(A).
\]
The algebra of all operators $A\in \A$ for which $PAP=A$ holds is
isomorphic to $M_3(\C)$. Let $r_1$ be the rank of $R_1$ and let $r_2$
be the rank of $R_2$. Clearly, we have $r_1+r_2=3$. The algebra of all
operators $B\in \B$ for which $B=R_1BR_1 +R_2BR_2$
is isomorphic to $M_{r_1}(\C) \oplus M_{r_2}(\C)$.
Therefore, $\phi$ induces a bijective transformation
\[
\psi : M_3(\C) \to M_{r_1}(\C) \oplus M_{r_2}(\C)
\]
which satisfies \eqref{E:triple3}.
We assert that either $r_1=3$ or $r_2=3$.
Suppose on the contrary that, for example, $r_1=2$ and $r_2=1$.
One can see that there are five rank-1 tripotents $P_1, \ldots ,P_5$
on the Hilbert space $\C^3$ such that $P_iP_jP_i=0$ $(i\neq j)$. Indeed,
choose an orthonormal basis $x,y,z$ in $\C^3$ and consider the operators
\[
(x+y)\ot x, \, y\ot (x+y), \, (1/2)(x-y)\ot (x-y), \,
z\ot (x+z), \, (x-z)\ot z.
\]
They fulfil the requirements. It follows that there are
five rank-1 tripotents in $M_2(\C)\oplus M_1(\C)$ with
similar properties. This readily implies that
there are four rank-1 tripotents $Q_1, \ldots, Q_4$ in $M_2(\C)$ for
which $Q_iQ_jQ_i=0$ $(i\neq j)$.
But this cannot happen. In fact,
applying a similarity transformation or the negative of a
similarity transformation
we can suppose that $Q_1=a\ot a$ for some unit vector $a$
in $\C^2$. Choose a unit vector $b\in \C^2$ which is orthogonal to $a$.
Since $Q_1Q_jQ_1=0$ $(j=2,3,4)$,
it follows that in the `vector-tensor-vector' representation of any of
$Q_2,Q_3,Q_4$ either the first or the second
component is a scalar multiple of $b$. Clearly, at least in
two of $Q_2,Q_3,Q_4$, the vector $b$ appears in the same component.
For example, suppose that $Q_2=c\ot b$ and $Q_3=d\ot b$.
Since $Q_2Q_3Q_2=0$, we obtain that either
$\la c,b\ra =0$ or $\la d,b\ra =0$. But this implies that either
$Q_2^2=0$ or $Q_3^2=0$ which contradicts the fact that $Q_2,Q_3$ are
nonzero tripotents.

We have proved that the induced transformation $\psi$ is a bijection
from $M_3(\C)$ onto itself which satisfies \eqref{E:triple3}.
Moreover, observe that we have also obtained that either $\psi(I)=I$
(this is the case if $r_1=3$) or $\psi(I)=-I$ (this is the case if
$r_2=3$). Without loss of generality we can suppose that $\psi(I)=I$.
It follows from \eqref{E:triple3} that
\[
\psi(A^2)=\psi(AIA)=\psi(A)I\psi(A)=\psi(A)^2
\]
which gives us that $\psi$ preserves the idempotents in both directions.
If $P,Q$ are idempotents in $M_3(\C)$ such that $PQ=QP=P$ (that is, if
$P\leq Q$), then we obtain
\begin{equation}\label{E:triple5}
\psi(P)\psi(Q)\psi(P)=\psi(P) \quad \text{and} \quad
\psi(Q)\psi(P)\psi(Q)=\psi(P).
\end{equation}
Since $\psi(P), \psi(Q)$ are idempotents, multiplying the second
equality in \eqref{E:triple5} by $\psi(Q)$ from the left and from the
right respectively, we find that
$\psi(Q)\psi(P)=\psi(P)\psi(Q)=\psi(P)$. So, we obtain that $\psi$
preserves the partial ordering $\leq$ between the idempotents in
$M_3(\C)$ in both directions.
We now apply a nice result of Ovchinnikov \cite{Ov}
describing the automorphisms of the poset of all idempotents on a
Hilbert space of dimension at least 3. It is a trivial corollary of his
result that our transformation $\psi$ is
orthoadditive on the set of all idempotents in $M_3(\C)$, that is, if
$P,Q$ are mutually orthogonal idempotents in $M_3(\C)$, then we have
$\psi(P+Q)=\psi(P)+\psi(Q)$.
Turning back to our original transformation $\phi$ we see that if
$P,Q\in \A$ are mutually orthogonal rank-1 idempotents, then
$\phi(P+Q)=\phi(P)+\phi(Q)$.

If $P$ is a rank-1 tripotent and $\lambda \in \C$ is a scalar, then we
have
\[
\phi(\lam P)=\phi(P (\lam P) P)=
\phi(P) \phi(\lam P) \phi(P)=h_P(\lam )\phi(P)
\]
for some scalar $h_P(\lam)\in \C$. This follows from the fact that
$\phi(P)$ has rank 1. We have
\[
h_P(\lam^2 \mu) \phi(P)=
\phi(\lam^2\mu P)=
\]
\[
\phi((\lam P)({\mu}P)(\lam P))=
\phi(\lam P)\phi({\mu}P)\phi(\lam P)=
h_P(\lam)^2h_P(\mu)\phi(P)
\]
which gives us that
\begin{equation}\label{E:triple7}
h_P(\lam^2 \mu)=h_P(\lam)^2h_P(\mu)
\end{equation}
for every $\lam, \mu \in \C$.
Choosing $\mu=1$, we see that $h_P(\lam^2)=h_P(\lam)^2$.
From \eqref{E:triple7} we now obtain
that $h_P$ is a multiplicative function.

We next assert that $h_P$ does not depend on $P$.
Let $Q\in \A$ be
a rank-1 tripotent with the property that $PQP\neq 0$.
We compute
\[
\phi((\lam P) (\mu^2 Q)(\lam P))=
\phi(\lam P) \phi(\mu^2 Q)\phi(\lam P)=
\]
\[
h_P(\lam)^2h_Q(\mu^2 )\phi(P)\phi(Q)\phi(P).
\]
On the other hand, we also have
\[
\phi((\lam P) (\mu^2 Q)(\lam P))=
\phi((\mu P) (\lam^2 Q)(\mu P))=
h_P(\mu)^2h_Q(\lam^2)\phi(P)\phi(Q)\phi(P).
\]
This yields that
\[
h_P(\lam)^2h_Q(\mu^2 )=
h_P(\mu)^2h_Q(\lam^2)
\]
and then we have $h_P=h_Q$.
If $PQP=0$, then we can choose a rank-1
tripotent $R\in \A$ such that $PRP\neq 0$ and $RQR\neq 0$.
Hence, we can infer that $h_P=h_R=h_Q$. This means that $h_P$ really
does not depend on $P$.
In what follows $h:\C \to \C$ denotes this common scalar function.

Let $A\in \A$ be arbitrary. Then we have
\[
\phi(P)\phi(\lam^2 A)\phi(P)=
\phi(P({\lambda}^2 A) P)=
\]
\[
\phi((\lam P) A (\lam P))=
\phi(\lam P)\phi(A)\phi(\lam P)=
h(\lam )^2\phi(P)\phi(A)\phi(P).
\]
Since this holds for every rank-1 tripotent $P$ on $X$ and
$\phi(P)$
runs through the whole set of rank-1 tripotents on $Y$, we
obtain that
$\phi(\lam^2 A)=h(\lam)^2\phi(A)$
for every $\lam \in \C$ which yields that
\[
\phi(\lam A)=h(\lam)\phi(A) \qquad (\lambda \in \C).
\]

We prove that $h$ is additive. Let $x,y\in X$ be linearly
independent vectors, and choose linear functionals $f,g\in X^*$ such
that $f(x)=1,f(y)=0$ and $g(x)=0, g(y)=1$.
Let $\lambda ,\mu \in \C$ be arbitrary and let $A=(\lam
x +\mu y)\ot (f+g)$, $P =x\ot f$, $Q=y\ot g$. By the orthoadditivity
property of $\phi$ we can compute
\[
h(\lam +\mu)\phi(A)=
\phi((\lam +\mu)A)=
\phi(A(P+Q)A)=
\]
\[
\phi(A)\phi(P+ Q)\phi(A)=
\phi(A)\phi(P)\phi(A)+ \phi(A)\phi(Q)\phi(A)=
\]
\[
\phi(AP A)+ \phi(A Q A)=
\phi(\lam A)+\phi(\mu A)=
(h(\lam)+h(\mu))\phi(A)
\]
and this proves that $h$ is additive.

We now verify that $\phi$ is additive.
Let $A,B\in \A$ be arbitrary and pick any rank-1 tripotent $P$
on $X$. We have $x\in X, f\in X^*$ such that $P=x\ot f$. We
compute
\begin{equation}\label{E:triple8}
\begin{gathered}
\phi(P)\phi(A+B)\phi(P)=
\phi(P(A+B) P)=
\phi(f((A+B)x) P)=
\\
h(f((A+B)x))\phi(P)=
h(f(Ax)) \phi(P)+h(f(Bx)) \phi(P)=
\\
\phi(f(Ax) P)+\phi(f(Bx) P)=
\phi(PA P)+\phi(P B P)=
\\
\phi(P)\phi(A) \phi(P)+\phi(P) \phi(B) \phi(P)=
\phi(P)(\phi(A)+\phi(B)) \phi(P).
\end{gathered}
\end{equation}
Since this holds true for every rank-1 tripotent $P$ on $X$, we
obtain that $\phi(A+B)=\phi(A)+\phi(B)$. Consequently, $\phi: \A \to \B$
is an additive bijection satisfying \eqref{E:triple3}.

Since every standard operator algebra $\mathcal R$ on a Banach space
is prime (this means that
for every $A,B \in \mathcal R$, the equality $A\mathcal R B=\{ 0\}$
implies $A=0$ or $B=0$), we can apply a result of Bre\v sar
\cite[Theorem 3.3]{Bresar} (also see \cite{Herstein2}) to obtain that
$\phi$ is necessarily a homomorphism, or an antihomomorphism,
or the negative of a homomorphism, or the negative of an
antihomomorphism.
In the homomorphic cases the satement follows from \cite{Semrl}, while
in the antihomomorphic cases one can apply analogous ideas (cf.
\cite[Proposition 3.1]{BresarSemrl}).
\end{proof}

\begin{remark}
We should explain why we have supposed in our theorem that
$\dim X\geq 3$. First, it
is easy to see that the conclusion does not hold true if $\dim
X=1$. Indeed, the function $z \longmapsto z|z|$ is a
multiplicative bijection of $\C$ which is not additive.
The place in the proof of our theorem where we have used that
$\dim X\geq 3$ is where we have applied Ovchinnikov's result.
If we knew that every bijective function $\phi: M_2(\C) \to M_2(\C)$
satisfying \eqref{E:triple3} preserves the mutual orthogonality between
rank-1 idempotents (that is, $\phi(P)\phi(Q)=\phi(Q)\phi(P)=0$ holds
whenever $P,Q\in M_2(\C)$ are rank-1 idempotents with $PQ=QP=0$),
then the use of this deep result could be avoided.
Unfortunately, we do not know whether $\phi$ has this preserving
property.
However, if we suppose that $\phi$ satisfies the stronger equality
\begin{equation}\label{E:triple12}
\phi(\{ ABC\})=\{ \phi(A)\phi(B)\phi(C)\}
\end{equation}
where $\{ABC\}$ denotes the so-called Jordan triple product $(1/2)
(ABC+CBA)$, then one can check that $\phi$ preserves
the mutual orthogonality between rank-1 idempotents and hence we
get the conclusion in Theorem also in the case $\dim
X=2$. Observe that if $\phi$ is additive, then the equation
\eqref{E:triple12} is equivalent to \eqref{E:triple3}.
\end{remark}

To conclude the paper we note that our approach in the considered
problem was mainly functional analytical. In our opinion, it is a
challenging
question how one can generalize the 'additive part' of our result for
general rings to obtain results similar to
Martindale's theorem.

% Bibliography
\bibliographystyle{amsplain}

\end{document}